\documentclass[reqno,10pt]{amsart}
\usepackage{amscd,amssymb,amsmath,color,mathtools}
\usepackage{paralist}
\usepackage{latexsym}
\usepackage[all]{xy}
\usepackage[colorlinks,allcolors=blue]{hyperref}
\usepackage[T1]{fontenc}
\usepackage{color}
\usepackage[english]{babel}
\usepackage{textcomp}
\usepackage{amsmath}
\usepackage{amssymb}
\usepackage{graphicx}
\usepackage{dsfont}

\newcommand{\lyxmathsym}[1]{\ifmmode\begingroup\def\b@ld{bold}
  \text{\ifx\math@version\b@ld\bfseries\fi#1}\endgroup\else#1\fi}

\def\..{{,\dots,}}
\def\:{{\colon}}

\def\int{{\rm int}}

\theoremstyle{definition}

\usepackage{booktabs} 

\graphicspath{{figures/}} 
 \usepackage{tikz} 
 \usetikzlibrary{calc}

\usepackage{tikz}
\usetikzlibrary{arrows}

\begin{document}

\author{Or Raz}

\date{\today}

\title{The geometry of ranked symplectic matroids}

\begin{abstract}
This paper is a continuation of my paper "Lattices of flats for symplectic matroids". We explore geometric constructions originating from the lattice of flats of ranked symplectic matroids. We observe that a ranked symplectic matroid always sits between two ordinary matroids and use this fact to prove that it has many of the same properties of ordinary matroids. We compute the dimension of its order complex using its Möbius function, We show that its matroid polytope is geometrically defined using its flats and connected to its Bergman fan. We finish by highlighting differences between its toric variety and the toric variety of an ordinary matroid, and give a partial proof of Mason's conjecture for ranked symplectic matroids.
\end{abstract}

\maketitle

\section{Introduction}

In \cite{[9]} Eva-Maria Feichtner and Bernd Sturmfels showed two important facts about the Bergman fan of an ordinary matroid.
\begin{enumerate}
    \item The order fan of its lattice of flats of an ordinary matroid refines its Bergman fan.
    \item The Bergman fan of an ordinary matroid is a subfan of the normal fan of its matroid polytope.
\end{enumerate}
These discoveries opened up the possibility of extending many important geometric results to other types of Coxeter matroids. This paper is an attempt at this for a special type of symplectic matroids.\\

In \cite{[2]} a ranked symplectic matroid was defined as the set of bases of a $C_n$ lattice, analogous to geometric lattices, as follows:
\subsection{Definition:}
\label{Definition1}
A collection $\mathcal{L}$ of subsets of $J=[n]\sqcup [n]^*$ is a $C_n$ lattice if:
\begin{enumerate}
\item $\emptyset,J\in\mathcal{L}$
\item $\forall A\in\mathcal{L}\setminus \left\{J\right\}$ we have $A\in P^{ad}(J)$.
\item $\forall A,B\in\mathcal{L}$ we have $A\cap B\in\mathcal{L}$.
\item For every $A\in\mathcal{L}$ let $\left\{ B_{1},...,B_{m}\right\} $
be the set of elements in $\mathcal{L}$ that cover $A$, then $\left\{ B_{1}\setminus A,...,B_{m}\setminus A\right\}$ is a partition of $J\setminus (A\sqcup A^{*})$.
\end{enumerate}
It was then shown that this family of lattices corresponds to a special type of symplectic matroids, called ranked symplectic matroids, for which they play the role of lattices of flats. \\

In contrast to what we did in \cite{[2]}, this paper defines ranked symplectic matroids using its enveloping matroid instead of a theoretic lattice of flats. This enables us to better understand the geometry of ranked symplectic matroids using the enveloping map, sending a pair of coordinates $x_i,x_{i^*}$ to $x_i-x_{i^*}$. Given a ranked symplectic matroid $\mathcal{S}$ and its enveloping matroid $\mathcal{M}$, the map $env$ not only enables us to find the matroid polytope of $\mathcal{S}$ in terms of that of $\mathcal{M}$ it also corresponds to a toric map between their Bergman fans.\\

The structure of this paper is as follows. In Section 2 we define a ranked symplectic matroid and present its beautiful structure. In Section 3 we study the Möbius function of a ranked symplectic matroid, showing that it alternates in sign using a recursive computation. In Section 4 we show that the matroid polytope of a ranked symplectic matroid can be defined by its lattice of flats, the same as ordinary matroids. In Section 5 we introduce the Bergman fan of a ranked symplectic matroid, proving that it is a subfan of the normal fan of its matroid polytope. Next, we study its group of Minkowski weights, characterizing the top-dimensional weights. We finish with Section 5 in which we prove Mason’s conjecture for rank $3$ ranked symplectic matroids.\\

\section{Preliminaries and basic properties}
In this paper, we make use of several properties and definitions concerning both ordinary and symplectic matroids. For background on ordinary matroids, we direct the reader to \cite{[1]}, and for symplectic matroids, to \cite{[8]}. We recall some definitions from \cite{[2]}. We say that a subset $A\subseteq J$ is admissible if $A\cap A^*=\emptyset$ and totally inadmissible if $A\triangle A^*=\emptyset$.

We denote by $I^{ad}(\mathcal{M})$ the family of independent sets of $\mathcal{M}$ with admissible closure. We call a symplectic matroid enveloped if it is the family of admissible bases of an ordinary matroid; our interest is in enveloped symplectic matroid with a specific type of enveloping matroid. One that induces a rank function on its enveloped symplectic matroid:

\subsection{Definition:}
\label{Definition2}
A loopless ordinary matroid $\mathcal{M}$ will be called admissible if its rank function is derived from its family of admissible independent sets $\mathcal{I}$ in the following way:

\[
r\left(A\right)=\max_{I\in\mathcal{I},\,I\subseteq A}\begin{cases}
\left|I\right|+2 & \exists\left\{ a,a^{*}\right\} \subseteq A\text{ s.t. }I\cup\left\{ a\right\} ,I\cup\left\{ a^{*}\right\} \in I^{ad}(\mathcal{M})\\
\left|I\right| & else
\end{cases}
\]
  
We now give a definition of ranked symplectic matroids, different from the definition found in \cite{[2]}, which is better suited to our needs.   
  
\subsection{Definition:}
\label{Definition3}
A ranked symplectic matroid $\mathcal{S}$ is the family of admissible bases of an admissible matroid $\mathcal{M}$. We always have that $\mathcal{M}$ is an enveloping matroid of $\mathcal{S}$.\\

We observe that for a ranked symplectic matroid $\mathcal{S}$, there always exists a unique minimal admissible enveloping matroid of the same rank. We will abuse the notation and use the term "enveloping matroid" to mean this minimal enveloping matroid.\\

We make use of the results presented in \cite{[2]} in Sections 4 and 5. To this end, we use remark [4.16] in \cite{[2]} to see Definition \ref{Definition3} is equivalent to the definition of ranked symplectic matroids in \cite{[2]}. Moreover, the lattice of flats $L(\mathcal{S})$ is equivalent to the $C_n$ lattice corresponding to $\mathcal{S}$.

\subsection*{Example:}
The simplest example of ranked symplectic matroids are uniform symplectic matroids, containing all admissible sets of size $k$ and denoted $U^*_{k,n}$.\\

The properties below were demonstrated in \cite{[2]} and likewise follow immediately from the definition of an admissible matroid.

\subsection{Lemma:}
\label{Lemma1}
\begin{enumerate}
\item If $\mathcal{M}$ is an admissible matroid then every flat is either admissible or totally inadmissible.
\item If $\mathcal{M}$ is an admissible matroid, then every atom of $\mathcal{M}$ is admissible.
\item If $\mathcal{M}$ is an admissible matroid, then $\mathcal{M}\setminus F$ is an admissible matroid for every inadmissible flat $F$ of $\mathcal{M}$.
\end{enumerate}

\subsection{Definition:}
\label{Definition4}
A strongly admissible set $A$ with respect to a matroid $\mathcal{M}$ is defined recursively. The empty set is strongly admissible, and $A$ is strongly admissible if there exists a strongly admissible set $B$ such that $A=B\cup \left\{a\right\}$ with $a,a^*\notin \overline{B}$. With $\overline{B}$ being the closure of $B$ in $\mathcal{M}$.

\subsection{Lemma:}
\label{Lemma2}
Let $\mathcal{M}$ be an admissible matroid. The family of strongly admissible sets of $\mathcal{M}$ is exactly $I^{ad}(\mathcal{M})$. 

\subsection*{proof:}
A strongly admissible set is independent, as every added element increases the rank of its closure.  Its closure is admissible because it cannot be totally admissible. On the other hand, $I^{ad}(\mathcal{M})$ is closed under taking subsets, implying that any $I\in I^{ad}(\mathcal{M})$ is strongly admissible.

The following Theorem shows that ranked symplectic matroids are indeed symplectic matroids.
\subsection{Theorem:}
\label{Theorem1}
If $\mathcal{M}$ is a matroid on $J$ in which every flat is admissible or totally inadmissible, then its admissible bases constitute a symplectic matroid.

\subsection*{proof:}
We use the definition of symplectic matroids by independent sets given in \cite{[3]}. 

\begin{enumerate}

\item Given a transversal $T\subset J$ we have $\mathcal{S}\setminus T= \mathcal{M}\setminus T$, and so an ordinary matroid.

\item Let $I_1,I_2\in I(\mathcal{S})$ be such that $\left|I_1\right|\lneq \left|I_2\right|$. As $I_1,I_2\in I(\mathcal{M})$ there exist $a\in I_2$ such that $ I_1\cup \left\{ a\right\}\in I(\mathcal{M})$. We assume that all such $a$ has $a^*\in I_1$, therefore $\overline{I_1}$ is not totally inadmissible and so admissible. Observe that $\left|I_1\setminus I_2^*\right|\lneq \left|I_2\setminus I_1^*\right|$ for the smaller independent sets $I_1\setminus I_2^*,I_2\setminus I_1^*$. Let $b\in I_2\setminus I_1^*$ be such that $(I_1\setminus I_2^*)\cup \left\{ b\right\}\in I(\mathcal{S})$. As $\overline{I_1}$ is admissible and $b,b^*\notin \overline{I_1}$ we have $I_1\cup \left\{ b^*\right\}\in I(\mathcal{S})$, satisfying the second part of the condition.

\end{enumerate}

Before diving into the main part of this section we introduce a very useful property which ranked symplectic matroids share with representable symplectic matroids.

\subsection{Definition:}
\label{Definition15}
A linear order $\omega$ on $J$ is called admissible if it satisfies the condition that for any two elements $a,b\in J$, the relation $a\le_\omega b$ holds if and only if $b^{*}\le_\omega a^{*}$.

\subsection{Lemma:}
\label{Lemma3}
Let $\mathcal{S}$ be a ranked symplectic matroid and $\omega$ be an admissible order on $J$ then the maximal base in the induced ordinary matroid $\mathcal{M}$ with respect to $\omega$ is admissible.

\subsection*{Proof:} 
Assume there is a maximal non-admissible base $B\in \mathcal{M}$ with respect to $\omega$ and let $a$ be the unique element with $a,a^*\in B$. As $B\setminus \left\{ a\right\}$ is strongly admissible, for any $b\in B\setminus \left\{ a,a^*\right\}$ we have  $(B\setminus \left\{ a\right\})\cup \left\{ b^*\right\}$ a base of $\mathcal{M}$. Assume that wlog $a^*\gneq a$, as $\omega$ is admissible, one of the following inequalities must hold:

\begin{enumerate}
\item $a\lneq b\lneq b^*\lneq a^*$
\item $b\lneq a\lneq a^*\lneq b^*$
\end{enumerate}

As $B\setminus \left\{ a\right\}$ is strongly admissible, we have in both cases that the two bases $B,(B\setminus \left\{ a,a^*\right\})\cup \left\{ b,b^*\right\}\in \mathcal{M}$ are incomparable.\\

Our main purpose of using ranked symplectic matroids is their correspondence with a lattice of flats and the operation of taking minors. Defined as follows:

\subsection{Definition:}
\label{Definition5}
Let $\mathcal{S}$ be a ranked symplectic matroid with an enveloping matroid $\mathcal{M}$. We defined the lattice of flats $L(\mathcal{S})$ of $\mathcal{S}$ to be the family of admissible flats of $\mathcal{M}$ plus $J$.

\subsection*{Example:}
In the uniform symplectic matroid $U_{k,n}^*$, the lattice of flats is formed by all admissible sets of size at most $k-1$, along with $J$.

\subsection{Definition:}
\label{Definition6}
Let $\mathcal{S}$ be a ranked symplectic matroid, $\mathcal{M}$ its enveloping matroid.
\begin{enumerate}
\item We define the deletion $\mathcal{S}\setminus a$ to be the ranked symplectic matroid corresponding to the admissible matroid $\mathcal{M}\setminus \left\{a,a^*\right\}$.
\item We define the contraction $\mathcal{S}/ a$ to be the ranked symplectic matroid corresponding to the admissible matroid $(\mathcal{M}/ a)\setminus \overline{a^*}$.
\end{enumerate}

 \subsection{Lemma:}
 \label{Lemma4}
If $\mathcal{M}$ is an admissible matroid of rank $d$, then the family of inadmissible flats of $\mathcal{M}$ forms a geometric lattice with the family of rank $d-1$ strongly admissible independent sets of $\mathcal{M}$ being its corresponding ordinary matroid.

\subsection*{proof:}
As every inadmissible flat is totally inadmissible, we have the family of inadmissible flats of $\mathcal{M}$ closed under intersection. As it is also upward closed in the lattice of flats of $\mathcal{M}$, it must be semi-modular and, therefore, a geometric lattice. To see that the corresponding ordinary matroid $\mathcal{M}'$ is the family of rank $d-1$ strongly admissible independent sets of $\mathcal{M}$, we first observe that it must be a sub-matroid of $\mathcal{M}$ of rank $d-1$. Every strongly independent set has an admissible closure in $\mathcal{M}$ and so its closure in $\mathcal{M}'$ will be the unique inadmissible set that covers it, making its rank equal to its size in $\mathcal{M}'$. In the other direction, if $B$ is a basis of $\mathcal{M}$, then it is admissible. We must have $rank_{\mathcal{M}'}(B)\geq rank_\mathcal{M}(B)$, so the closure of $B$ in $\mathcal{M}$ is admissible, making $B$ strongly admissible.
\\


We now observe that every admissible basis of an admissible matroid is "almost" strongly admissible, meaning it is obtained from a strongly admissible independent set by adding an element. Another use for the family of rank $d-1$ strongly admissible set, using this observation, is the following "almost" matroid isomorphism of the restrictions by inadmissible pairs. 

\subsection{Lemma:}
\label{Lemma5}
Let $\mathcal{S}$ be a ranked symplectic matroid of rank $\geq 3$ with an enveloping matroid $\mathcal{M}$, then for every $a\in J$ the following map is an isomorphism of sets, $\psi :T(\mathcal{S}/ a)\rightarrow T(\mathcal{S}/ a^*)$

\[
\psi\left(A\cup\left\{ b\right\} \right)=\begin{cases}
A\cup\left\{ b\right\}  & \textbf{if } A\cup\left\{ b\right\} \text{ is strongly admissible}\\
A\cup\left\{ b^{*}\right\}  & else
\end{cases}
\]
 
With $A\cup\left\{ b\right\}$ being the decomposition of a basis to a maximal strongly admissible independent set and an element.
\subsection*{proof:}
Using Lemma \ref{Lemma4} we have that the matroid $\mathcal{M}'$ consists of strongly admissible independent sets of rank $d-1$ having only totally inadmissible flats. Therefore, the identity map on $J\setminus \left\{a,a^*\right\}$ induces an isomorphism on the set of strongly admissible bases of $T(\mathcal{S}/ a)$ and $T(\mathcal{S}/ a^*)$. As the decomposition $A\cup\left\{ b\right\}$ is unique, it is enough to show that $\psi$ is surjective. Let $A\cup\left\{ b\right\}$ be an admissible, but not strongly admissible basis of $T(\mathcal{M}/ a^*)$. We have $\overline{A\cup\left\{ a\right\}},\overline{A\cup\left\{ a^*\right\}}$ non-maximal flats that cover $\overline{A}$ in $\mathcal{M}$. Consequently, $\overline{A\cup\left\{ a\right\}}=\overline{A}\cup(\overline{A}\setminus \overline{A\cup\left\{ a^*\right\}})^* $ and if $b\in \overline{A\cup\left\{ a^*\right\}}$ then $b^*\in \overline{A\cup\left\{ a\right\}}$.

\subsection*{Remark:}
$T((\mathcal{M}/ a))\setminus a^*)$ and $T((\mathcal{M}/ a^*)\setminus a)$ are never isomorphic in the sense of matroids. The isomorphism will never be induced by a bijection on the ground set.

\subsection{Lemma:}
\label{Lemma6}
Let $F'\neq J$ be an inadmissible flat of an admissible matroid $\mathcal{M}$ and $\mathcal{F}$ the set of admissible flats covered by $F'$.
\begin{enumerate}
\item Every admissible flat $F\in \mathcal{F}$ is a transversal in $\mathcal{M}\setminus F$, i.e., $F'=F\sqcup F^*$.
\item $[\overline{\left\{a,a^*\right\}},F']\cong [\overline{\left\{a\right\}},F]$ for every $a\in F\in \mathcal{F}$.
\end{enumerate}

\subsection*{proof:} 
To see (1) assume $F'\in \mathcal{F}$ and let $I$ be a strongly admissible independent set that spans $F$. If $a,a^*\notin F'$, then $I\cup \left\{a\right\}$ is strongly admissible and thus has an admissible closure. As $F$ is inadmissible, we must have $a,a^*\notin F$.

We define the isomorphism in (2) as the map that sends each element in $[\overline{\left\{i\right\}},F]$ to the unique inadmissible element that covers it. This map is surjective as every inadmissible flat containing $a$ covers some admissible element containing $a$. It is also injective because for each $G\in [\overline{\left\{a,a^*\right\}},F']$ there is at most one transversal of $\mathcal{M}\setminus G$ contained in $[\overline{\left\{a\right\}},F]$.\\

We finish this section with a discussion on connectivity. Given an ordinary matroid $\mathcal{M}$ we say that $\mathcal{M}$ is connected if it cannot be written as a direct sum of ordinary matroids on disjoint groundsets. For ordinary matroids, we see that two elements of the ground set are in the same connected component iff they are both contained in a circuit together. For an enveloping matroid of a ranked symplectic matroid, we have the following:

\subsection{Lemma:}
\label{Lemma7}
The enveloping matroid $\mathcal{M}$ of a ranked symplectic matroid on the groundset $J$ is either connected or of rank $2$.

\subsection*{proof:}
It is enough to show that if $rank(\mathcal{M})\geq 3$ then all pairs of elements are contained in a circuit together (we can assume that $rank(\mathcal{M})\geq 2$ as all rank $1$ matroids are connected). Let $F$ be the rank $1$ flat of $\mathcal{M}$ and $a\in F$. As $rank(\mathcal{M})\geq 3$, we find that $F$ is not a transversal and that $F^*$ is also a flat. Therefore, the set $\left\{a,a^*,b,b^*\right\}$ is a circuit of $\mathcal{M}$ for every $b\notin F\cup F^*$. As two elements of $J$ contained in the same rank $1$ flat constitute a circuit we are done.\\

To finish the discussion, we note that there is no good way to define connectedness for ranked symplectic matroids, as being a ranked symplectic matroid is not closed under taking direct sum (and in the same way their lattices are not closed under taking direct product). Take, for example, the two ranked symplectic matroids $\mathcal{S}_1=\left\{\left\{1,2^*\right\},\left\{1^*,2\right\},\left\{1^*,2^*\right\}\right\}$ and $\mathcal{S}_2=\left\{\left\{3\right\},\left\{3^*\right\}\right\}$. We do not have $\mathcal{S}_1\oplus \mathcal{S}_2$ a ranked symplectic matroid. Moreover, we can see from Lemma \ref{Lemma5} that the direct sum of enveloping matroids is never an enveloping matroid.

\section{The Möbius function of a ranked symplectic matroid}

One of the main tools of understanding the structure of a poset is its Möbius function; it is even more important in the case of shellable posets, as it is then equal to the dimension of the homology of its order complex. We refer the reader to \cite{[16]} for further details on the theory of shellable posets.\\

We begin by introducing the Möbius function of a finite poset, together with the properties required for our discussion. Proofs and additional details can be found in \cite{[17]}.

\subsection{Definition:}
\label{Definition7}
Let $P$ be a finite poset. The Möbius function $\mu:P\rightarrow \mathbb{Z}$ is defined recursively by
\[
\mu\left(x,y\right)=\begin{cases}
1 & y=x\\
-\sum_{x\leq z\lneq y}\mu\left(x,z\right) & y\geq x\\
0 & \text{else}
\end{cases}
\]

We will need the following general Theorems on the Möbius function of a finite lattice and the specific Theorem on the Möbius function of a geometric lattice:

 \subsection{Theorem (Weisner):}
 \label{Theorem2}
 Let $L$ be a finite lattice and $a\in A(L)$, an element covering $\hat{0}$. Then:
 
 \[
\sum_{F\in L:F\vee a=1_{L}}\mu\left(0_{L},F\right)=0
\]

\subsection{Theorem (Boolean Expansion Formula) :}
\label{Theorem3}
 Let $L$ be a finite lattice. Then:
 
 \[
\mu (L)=\sum_{\begin{matrix}B\in2^{A(L)}\\
\vee B=1_{L}
\end{matrix}}(-1)^{\left|B\right|}
\]

 \subsection{Theorem:}
 \label{Theorem4}
 Let $\mathcal{M}$ be an ordinary loopless matroid on the ground set $E$, then the following hold:
 \begin{enumerate}
\item $\mu(\mathcal{M})= \mu(\mathcal{M}\setminus a)-\mu(\mathcal{M}/a)$ for all $a\in E$ not a coloop.
\item $\mu(\mathcal{M})=-\underset{\begin{smallmatrix}a\notin F\in L\\
rank\left(F\right)=d-1
\end{smallmatrix}}{\sum}\mu\left(\hat{0},F\right)$
\end{enumerate}



Using the Boolean Expansion Formula we can formulate the difference between $\mu (\mathcal{M})$ and $\mu (\mathcal{S})$. Let $X\subsetneq 2^{A(L(\mathcal{M}))}$ denote the set of inadmissible non-spanning sets of atoms of the enveloping matroid $\mathcal{M}$, then

\[
\mu\left(\mathcal{S}\right)=\mu\left(\mathcal{M}\right)+\sum_{B\in X} (-1)^{\left|B\right|}
\]

 We can use this to prove the following Lemma.

\subsection{Lemma:}
\label{Lemma8}
Let $\mathcal{S}$ be a ranked symplectic matroid of rank $d$ and let $\mathcal{F}_a(k)$ be the set of rank $k$ flats of $\mathcal{S}$ containing $a$. Then for every $a\in J$ we have:

\[
\sum_{F\in\mathcal{F}_{a}\left(k\right)}\mu\left(F\right)=\sum_{F\in\mathcal{F}_{a^{*}}\left(k\right)}\mu\left(F\right)
\]

\subsection*{proof:} 
It enough to prove this for $k=d-1$ as every flat is of that rank in some truncation of $\mathcal{S}$. We use Theorem \ref{Theorem4} to obtain the following for the enveloping matroid $\mathcal{M}$:

\[
\underset{\begin{smallmatrix}a^*\notin F\\
r\left(F\right)=d-1
\end{smallmatrix}}{\sum}\mu\left(\hat{0},F\right) -\underset{\begin{smallmatrix}a\notin F\\
r\left(F\right)=d-1
\end{smallmatrix}}{\sum}\mu\left(\hat{0},F\right)=-\mu(\mathcal{M})+\mu(\mathcal{M})=0
\]

Observe that every flat not containing both $a$ and $a^*$ appears in both sums, and every flat containing both $a$ and $a^*$ does not appear in both sums. Therefore, we are left with only flats containing $a$ or $a^*$, but not both, which are always admissible and so flats of $\mathcal{S}$.\\

The Möbius function of a ranked symplectic matroid also has a deletion-contraction rule, similarly to ordinary matroids, we call $\left\{a,a^*\right\}$ a coloop of $\mathcal{S}$ if $\mathcal{S}\setminus \overline{\left\{a,a^*\right\}}=U^*_{m,m}$ for some $m\in \mathbb{N}$.
 
\subsection{Theorem:}
\label{Theorem5}
Let $\mathcal{S}$ be a ranked symplectic matroid of rank $\geq 3$, if $\left\{a,a^*\right\}$ is not a coloop of $\mathcal{S}$ then:
\[
\mu(\mathcal{S})= \mu(\mathcal{S}\setminus \overline{\left\{a,a^*\right\}})-\mu(\mathcal{S}/\overline{a})-\mu(\mathcal{S}/\overline{a^*})
\]

\subsection*{proof:} 
As deleting all elements parallel to $a$ or $a^*$ is a lattice isomorphism, we will prove the Theorem for case $\overline{a}=\left\{a\right\}$. We first use Theorem \ref{Theorem4} twice with $a$ and then $a^*$, observing that by Definition \ref{Definition6} the matroids $\mathcal{M},\mathcal{M}\setminus a$ and $\mathcal{M}/\overline{a}$ are always without coloops. 

\[
\mu(\mathcal{S})-\sum_{B\in X} (-1)^{\left|B\right|}= \mu(\mathcal{M})=\mu(\mathcal{M}\setminus a)-\mu(\mathcal{M}/a)
\]
\[
= \mu(\mathcal{M}\setminus \left\{a,a^*\right\})-\mu((\mathcal{M}\setminus a)/a^*)-\mu((\mathcal{M}/a)\setminus a^*)+\mu(\mathcal{M}/\left\{a,a^*\right\})
\]

 To complete the proof, we decompose $X$ into a sum of four subsets. Let $X_{-a-a^*}$ be the sets that contain neither $a$ or $a^*$, $X_{a-a^*}$ the sets that contain $a$ but not $a^*$, $X_{a^*-a}$ the sets that contain $a^*$ but not $a$ and $X_{a+a^*}$ the sets that contain both $a$ and $a^*$. We obtain the following four equalities, finishing the proof.
 
\begin{enumerate}
\item 
$\mu\left(\mathcal{S}\setminus \left\{a,a^*\right\}\right)=\mu(\mathcal{M}\setminus \left\{a,a^*\right\})+\sum_{B\in X_{-a-a^*}} (-1)^{\left|B\right|}$
\item 
$\mu\left(\mathcal{S}/a\right)=\mu((\mathcal{M}/a)\setminus a^*)+\sum_{B\in X_{a-a^*}} (-1)^{\left|B\right|}$
\item 
$\mu\left(\mathcal{S}/a^*\right)=\mu((\mathcal{M}\setminus a)/a^*)+\sum_{B\in X_{a^*-a}} (-1)^{\left|B\right|}$
\item 
$\mu(\mathcal{M}/\left\{a,a^*\right\})+\sum_{B\in X_{a+a^*}} (-1)^{\left|B\right|}=\sum_{B\in 2^{J\setminus\left\{a,a^*\right\}}}  (-1)^{\left|B\right|}=0$
\end{enumerate}

\subsection{Theorem:}
\label{Theorem6}
The Möbius function of a ranked symplectic matroid $\mathcal{S}$ is nonzero and alternates in sign, precisely, for every $F_1\leq F_2\in L(\mathcal{S})$:
\[
(-1)^{rank(F_1)-rank(F_2)}\mu_{L(\mathcal{S})}(F_1,F_2)>0
\]

\subsection*{proof:} 

For any $F_2\neq J$ the Theorem follows from the same property of the enveloping matroid and by taking minors it is enough to prove the Theorem for the case $\mu_{L(\mathcal{S})}(F_1,F_2)=\mu(\mathcal{S})$. Using Theorem \ref{Theorem5} we are left to prove that the Theorem holds for $\mathcal{S}$ with every element a coloops. In that case, every atom of $L(\mathcal{S})$ is a singleton and therefore $(\mathcal{M}\setminus a)/a^*$ is an admissible matroid. The following holds:

\[
\mu(\mathcal{S})-\sum_{B\in X} (-1)^{\left|B\right|}= \mu(\mathcal{M}\setminus \left\{a,a^*\right\})-\mu((\mathcal{M}\setminus a)/a^*)-\mu((\mathcal{M}/a)\setminus a^*)+\mu(\mathcal{M}/\left\{a,a^*\right\})
\]

\[
\Rightarrow \mu(\mathcal{S})= \mu(\mathcal{M}\setminus \left\{a,a^*\right\})-\mu\left(\mathcal{S}/a\right)-\mu\left(\mathcal{S}/a^*\right)- \mu(T(\mathcal{M}\setminus \left\{a,a^*\right\}))+ \mu(T(\mathcal{S}\setminus \left\{a,a^*\right\}))
\]

Most summands are computed in the same way as in Theorem \ref{Theorem5}, notice that as $\left\{a,a^*\right\})$ is a coloop, we have $T(\mathcal{S}\setminus \left\{a,a^*\right\})=U^*_{n-1,n-1}$ and therefore:

\[
- \mu(T(\mathcal{M}\setminus \left\{a,a^*\right\}))+ \mu(T(\mathcal{S}\setminus \left\{a,a^*\right\}))= -\sum_{k=0}^{n-4}\left(-1\right)^{k}n{2n-4 \choose k}=\left(-1\right)^{n}\frac{n}{2}{2n-4 \choose n-2}
\]

We then have the sign of $\mu(\mathcal{S})$ correct and finish by induction, for which the base case of $rank(\mathcal{S})=2$ is obvious. 

\section{Matroid polytopes}
Another cryptomorphic definition of matroids, one that also holds for symplectic matroids, is the matroid polytope definition given by the Gelfand-Serganova Theorem. 

We review the connection between the matroid polytopes of a ranked symplectic matroid and the matroid polytope of its enveloping matroid. Let $A$ be an $m$-element subset of the ground set $J$. We represent $\{ a_1,\ldots ,a_m\} =A\subset [n]$ in $\mathbb{R^{\mathnormal{2n}}}$ and $\mathbb{R^{\mathnormal{n}}}$, respectively, by the following two corresponding sums of $m$ unit vectors:

\[
e_A=\sum_{a_i\in A}e_{a_i}\in \mathbb{R}^{2n}
\]
\[
e^{\pm}_A=\sum_{a_i\in A\cap [n]} e_{a_i}-\sum_{a_i\in A\cap[n]^*} e_{a_i^*} \in \mathbb{R}^{n}
\]

\subsection{Definition:} 
\label{Definition8}
The matroid polytope $P(\mathcal{M})$ of $\mathcal{M}$ a matroid/symplectic matroid is the convex hull of $\left\{e_{B}\mid B\in\mathcal{M}\right\}$ and $\left\{e^{\pm}_{B}\mid B\in\mathcal{M}\right\}$ respectively.\\

For completeness, we present a special case of the Gelfand–Serganova Theorem for symplectic matroids. Further details can be found in \cite{[8]}.

\subsection{Theorem:}
\label{Theorem7}
Let $\mathcal{B}$ be a set of admissible $k$ sets in $J$. Let $P(\mathcal{B})$ be the convex hull of
the points $e^{\pm}_A$ with $A\in \mathcal{B}$.
Then $e^{\pm}_A$ are the vertices of $P(\mathcal{B})$ for all $A\in \mathcal{B}$. Moreover, the set $\mathcal{B}$ is the collection of bases
of a symplectic matroid on $J$ if and only if all edges (that is, one-dimensional faces) of $P(\mathcal{B})$ are
parallel to a vector of the form $e^{\pm}_i-e^{\pm}_j$ for $i,j\in J$.\\

To establish a connection between the matroid polytope of a ranked symplectic matroid and the matroid polytope, we use the following enveloping map:
\[
env:\mathbb{R}^{2n}\rightarrow \mathbb{R}^{n}
\]
\[
env(x_1,...,x_n,x_{1^*},...,x_{n^*})=(x_1-x_{1^*},...,x_n-x_{n^*})
\]
We would have liked $env(P(\mathcal{M}))=P(\mathcal{S})$, for a ranked symplectic matroid $\mathcal{S}$ and its enveloping matroid $\mathcal{M}$, but unfortunately that is not the case. We observe the following connection:

\subsection{Theorem:}
\label{Theorem8}
Let $\mathcal{S}$ be a ranked symplectic matroid and $\mathcal{M}$ its enveloping matroid. The image of an inadmissible basis $B$ of $\mathcal{M}$ under the enveloping map is in $P(\mathcal{S})$ iff the closure of $B\setminus \left\{a,a^*\right\}$ is not a transversal in $J\setminus \left\{a,a^*\right\}$. With $a,a^*$ being the inadmissible pair contained in $B$. 

\subsection*{Proof:}
Let $B\in \mathcal{M}$ be an inadmissible base, we have to prove:
\[
env(e_B)=e_{B\setminus (B\triangle B^*)}=e_{B\setminus \left\{ a,a^*\right\}} \in P(\mathcal{S})  
\]

As $B$ is inadmissible, we have $B\setminus \left\{a,a^*\right\}$ a not maximal, strongly admissible independent set and therefore it has an admissible closure $F$. Using our discussion of minors of ranked symplectic matroids, this is equivalent to proving that $\mathcal{S}/F$ is of rank $2$ and the origin is contained in $P(\mathcal{S}/F)$. As $F$ is not a transversal in $J\setminus \left\{a,a^*\right\}$, we have $\mathcal{S}/F$ a ranked symplectic matroid of rank $2$. The proper flats of a rank $2$ ranked symplectic matroids are admissible, and therefore there must exist $b,c\in \mathcal{S}/F$ such that $\left\{ b,c^*\right\},\left\{ c,b^*\right\}$ are bases of $\mathcal{S}/F$ and therefore the origin is contained in their convex hull.

On the other hand, if $B\setminus \left\{a,a^*\right\}$ is a transversal in $J\setminus \left\{a,a^*\right\}$, then $\mathcal{S}/F=\left\{a,a^*\right\}$ which is of rank $1$ and the image of $B$ is not in $P(\mathcal{S})$. In particular, we see that $env(P(\mathcal{M}))=P(\mathcal{S})$ iff $L(\mathcal{M})$ does not contain any flat of the form $J\setminus \left\{a,a^*\right\}$ for some $a\in [n]$. 

The second goal of this section is to produce a description of the matroid polytope of a ranked symplectic matroid in terms of its flats in a way similar to the one that is known for ordinary matroids. We start with some definitions and the ordinary case.

In the case of a rank $r$ ordinary matroid on $2n$ elements the matroid polytope is a subset of the simplex:
\[
\Delta =\left\{(x_1,...,x_{2n})\in \mathbb{R}^{2n} \mid x_1\geq 0,...,x_{2n}\geq 0,\sum_{i=1}^{2n} x_i=r\right\}
\]
In our case, the matroid polytope of a rank $r$ symplectic matroid on $2n$ elements is a subset of the cross polytope:
\[
\diamondsuit =\left\{(x_1,...,x_n)\in \mathbb{R}^n \mid \sum_{i=1}^{n}\left|x_{i}\right|\leq r\right\}
\]

We have the following characterization of its matroid polytope using its lattice of flats:

\subsection{Theorem:}
\label{Theorem9}
The matroid polytope of a rank $r$ ordinary matroid $\mathcal{M}$ on the ground set $[n]$ is equal to the following subset of the $n$-simplex:
\[
P(\mathcal{M})=\left\{(x_1,...,x_n)\in \Delta \mid \sum_{i\in F}x_i\leq rank_\mathcal{M}(F),\forall F\in L(\mathcal{M})\right\}
\]
Furthermore, if $\mathcal{M}$ is connected, then the inequalities defining the facets of $P(\mathcal{M})$ correspond to flats such that both $\mathcal{M}\setminus F$ and $\mathcal{M}/ F$ are connected.

\subsection{Theorem:}
\label{Theorem10}
The matroid polytope of a rank $r$ ranked symplectic matroid $\mathcal{S}$ on the ground set $J$ is equal to the following subset of the $n$-cross polytope:
\[
P(\mathcal{S})=\left\{(x_1,...,x_n)\in \diamondsuit \mid \sum_{\left\{i,i^*\right\}\cap F\neq \emptyset}sgn(i)\cdot x_i\leq rank_\mathcal{S}(F)-\phi (F),\forall F\in L(\mathcal{S})\right\}
\]
With $\phi (F)=\begin{cases}
1 &  \textbf{if }|F|=n-1 \textbf{ if } rank(F)=r-2\\
0 & else\\
\end{cases}$
\subsection*{Proof:}
We first observe that if $\mathcal{M}$ is the enveloping matroid of $\mathcal{S}$ then the admissible flats of $L(\mathcal{M})$ are exactly the flats in $L(\mathcal{S})\setminus \left\{J\right\}$ and we have: 
\[
env(P(\mathcal{M}))=\left\{(x_1,...,x_n)\in \diamondsuit \mid \sum_{\left\{i,i^*\right\}\cap F\neq \emptyset}sgn(i)\cdot x_i\leq rank_{\mathcal{S}}(F),\forall F\in L(\mathcal{S})\right\}
\]
We are now left with eliminating the contribution of the non admissible bases of $\mathcal{M}$. By Theorem \ref{Theorem7} we have every inadmissible basis $B$ of $\mathcal{M}$ with $env(e_B)\notin P(\mathcal{S})$ admitting the equality $\sum_{\left\{i,i^*\right\}\cap F\neq \emptyset}sgn(i)\cdot x_i= rank_\mathcal{S}(F)$ for some flat $F$ such that $\phi (F)=1$ and therefore not admitting our amended list of inequalities. Furthermore, let $B$ be any basis for $\mathcal{S}$ maximizing $\sum_{\left\{i,i^*\right\}\cap F\neq \emptyset}sgn(i)\cdot x_i$ for a flat $F$ such that $\phi (F)=1$. We must have a subset $I\subset B$ that spans $F$ and since $F$ is a transversal in some $J\setminus \left\{a,a^*\right\}$ we must have $B=I\cup \left\{a,b\right\}$ or $B=I\cup \left\{a^*,b\right\}$ for some $b\in J$ with $b^*\in F$. Let $(x_1,...,x_n)=e^{\pm}_B$:
\[
\sum_{\left\{i,i^*\right\}\cap F\neq \emptyset}sgn(i)\cdot x_i=rank_\mathcal{S}(F)+sgn(b^*)\cdot x_b=rank_\mathcal{S}(F)-\phi (F)
\]
And we therefore do not eliminate any admissible bases.

We finish the proof by showing that $\sum_{\left\{i,i^*\right\}\cap F\neq \emptyset}sgn(i)\cdot x_i=rank_\mathcal{S}(F)-\phi (F)$ is a facet defining equality in $P(\mathcal{S})$. Using the same argument as before, we see that the hyperplane defined by the bases maximizing $F$ is the same as the bases maximizing the totally inadmissible flat $F'=J\setminus \left\{a,a^*\right\}$ covering $F$ in $L(\mathcal{M})$. To see $F'$ is a facet that defines the flat of $P(\mathcal{M})$, we use Lemma \ref{Lemma6}. If $rank(F')\geq 2$, then $\mathcal{M}\setminus F'$ is an admissible matroid and so connected. The ordinary matroid $\mathcal{M}/ F'$ is equal to a rank $1$ matroid with bases $\left\{a\right\},\left\{a^*\right\}$, which is connected. As the hyperplane defined by $F$ in $P(\mathcal{S})$ is the image under $env$ of a facet that defines the hyperplane of $P(\mathcal{M})$, we are done. If $rank(F')= 2$ then the proper flats of $\mathcal{S}$ are exactly:
\[
F,F^*,\left\{a\right\},\left\{a^*\right\},F\cup \left\{a\right\},F\cup \left\{a^*\right\},F^*\cup \left\{a\right\},F^*\cup \left\{a^*\right\}\subseteq \mathbb{R}^n
\]
And it is easily checked by hand that $\sum_{\left\{i,i^*\right\}\cap F\neq \emptyset}sgn(i)\cdot x_i= rank_\mathcal{S}(F)-1$ is a facet defining equality.\\

We finish this section with a discussion on the dimension of the matroid polytope of a ranked symplectic matroid. For ordinary matroids, the following Theorem holds:
\subsection{Theorem:}
\label{Theorem11}
The dimension of the matroid polytope of an ordinary matroid $\mathcal{M}$ defined on the ground set $[n]$ is $dimP(\mathcal{M})=n-c(\mathcal{M})$, with $c(\mathcal{M})$ being the number of connected components of $\mathcal{M}$.\\

In the case of a ranked symplectic matroid $\mathcal{S}$, we do not have a good definition of connected components. However, we will see that in the same way as its enveloping matroid, every ranked symplectic matroid of rank $\geq 3$ corresponds to a full-dimensional matroid polytope. 

\subsection{Theorem:}
\label{Theorem12}
Let $\mathcal{S}$ be a ranked symplectic matroid on the ground set $J$. If $\mathcal{S}$ is not isomorphic to $\left\{\left\{1,2\right\},\left\{1^*,2^*\right\}\right\}$, then $dimP(\mathcal{S})=n$. 

\subsection*{Proof:}
Let $rank(\mathcal{S})\geq 3$ and $i\in J$. As $\left\{i,i^*\right\}$ is independent in the enveloping matroid, it contains a basis $\left\{i,i^*\right\} \subseteq B$. We can therefore extend $B\setminus \left\{i\right\}$ to a basis $(B\setminus \left\{i\right\})\cup \left\{j\right\}$ of $\mathcal{S}$. As $rank(\mathcal{S})\geq 3$ we have $i,i^*\notin \overline{(B\setminus \left\{i,i^*\right\})\cup \left\{j\right\}}$ and therefore both $(B\setminus \left\{i\right\})\cup \left\{j\right\}$ and $(B\setminus \left\{i^*\right\})\cup \left\{j\right\}$  are bases of $S$. We obtain a segment in $P(\mathcal{S})$ parallel to $e_i$ for all $i\in [n]$.

If $rank(\mathcal{S})=2$ and $\mathcal{S}\ncong \left\{\left\{1,2\right\},\left\{1^*,2^*\right\}\right\}$ then we can still find such a $j$ if $\left|L(\mathcal{S})\right|\neq 4$.

If $rank(\mathcal{S})=2$, $\mathcal{S}\ncong \left\{\left\{1,2\right\},\left\{1^*,2^*\right\}\right\}$, and $\left|L(\mathcal{S})\right|= 4$ then we must have $\left|J\right|\geq 6$. We then have $\left\{i,j\right\},\left\{i^*,j^*\right\}\in \mathcal{S}$ for every $i,j\in J$ that corresponds to a vector set in $P(\mathcal{S})$ generating $\mathbb{R}^n$.

\section{The Bergman fan of a ranked symplectic matroid}
One of the most extraordinary developments in matroid theory in recent years is the study of its Bergman fan, its associated toric variety, and its Chow ring. We follow the definitions and constructions given in \cite{[5]}, to keep this work self-contained we will repeat the necessary results.

\subsection{Definition:}
\label{Definition9}
Let $\mathcal{M}$ be an ordinary matroid on the ground set $[n]$. Its Bergman fan is the following collection of cones in $\mathbb{P}^{n}=(\mathbb{Z}^n/\left\langle e_{[n]}\right\rangle)\otimes \mathbb{R}$:
\[
\left\{[Span_{\mathbb{R}_+}(\mathcal{F})]\mid \mathcal{F}=\left\{F_1\subsetneq...\subsetneq F_k\right\}\subset L(\mathcal{M})\setminus \left\{\emptyset,J\right\} \right\} 
\]
with $Span_{\mathbb{R}_+}(\mathcal{F})=\left\{\sum_{F\in\mathcal{F}}a_{F}\cdot e_{F}\mid a_{F}\in\mathbb{R}_{\geq0}\right\}$ and $[Span_{\mathbb{R}_+}(\mathcal{F})]$ its image under the quotient map $\mathbb{R}^n\rightarrow \mathbb{P}^{n}$. In some Theorems we will need this realization of the order fan to be in $\mathbb{R}^{n}$ instead of $\mathbb{P}^n$, in which case we will denote it $B(\mathcal{M})_{\mathbb{R}}$ instead, omitting the quotient map from the definition.

Some of our readers may know the Bergman fan of $\mathcal{B}$ to be its tropical linear space. This definition was shown in \cite{[13]} to be its fine subdivision.

We define the Bergman fan of the ranked symplectic matroid $\mathcal{S}$, in a similar way, to be the order fan of the proper part of its lattice of flats with the realization discussed in the beginning of Section 5.

\subsection{Definition:}
\label{Definition10}
The Bergman fan $B(\mathcal{S})$ of the ranked symplectic matroid $\mathcal{S}$ on $J$ is a fan in $\mathbb{R}^n$ consisting of the following set of cones in $\mathbb{R}^n$:
\[
\left\{Span_{\mathbb{R}_+}(\mathcal{F})\mid \mathcal{F}=\left\{F_1\subsetneq...\subsetneq F_k\right\}\subset L(\mathcal{S})\setminus \left\{\emptyset,J\right\} \right\} 
\]
With $Span_{\mathbb{R}_+}(\mathcal{F})=\left\{\sum_{F\in\mathcal{F}}a_{F}\cdot e_{F}^{\pm}\mid a_{F}\in\mathbb{R}_{\geq0}\right\}$.\\

\subsection{Remark:}
A fan is defined as a collection of polyhedral cones with apex at the origin such that every face of a cone belongs to the fan, and the intersection of any two cones is a face of both. The Bergman fans $B(\mathcal{M})$ and $B(\mathcal{S})$ are fans because of their lattice-of-flats definition.\\

In \cite{[2]} it was proven that the lattice of flats of a ranked symplectic matroid $L(\mathcal{S})$ is shellable and therefore $B(\mathcal{S})$ it is connected in co-dimension $1$. Another property we will need is unimodularity: 

\subsection{Definition:}
\label{Definition11}
A fan in $\mathbb{R}^n$ will be called unimodular if, for each cone of the fan, the first lattice point along each of the rays forms a basis for the lattice in the linear span of the cone.\\

The fan $B(\mathcal{S})$ is unimodular because its cones are generated by subsets of the standard basis vectors of $\mathbb{R^{\mathnormal{n}}}$. This means that each cone corresponds to a coordinate subspace, and the generating vectors form part of a unimodular matrix (determinant $\pm 1$), ensuring that the fan is unimodular.\\

Next, we show that there is a connection between the Bergman fan and the matroid polytope of a ranked symplectic matroid. We first present the necessary results from \cite{[9]}. Let $\mathcal{M}$ be an ordinary matroid on the ground set $[n]$ and $\nu \in \mathbb{R}^n$. Using the Gelfand-Serganova Theorem, we can define $\mathcal{M}_{\nu}$ as the ordinary matroid corresponding to the face of $P(\mathcal{M})$ on which the linear functional $\sum^n_{i=1}x_i\nu_i$ attains its maximum. The fibers of the map $\nu→ M_{\nu}$ are open polyhedral cones, whose closures form a complete fan. In \cite{[9]} it was shown that $B(\mathcal{M})_{\mathbb{R}}$ is a refinement of a sub-fan of the normal fan of $P(\mathcal{M})$, and more specifically, the following Theorem holds:

\subsection{Theorem (Proposition 2.5. \cite{[9]}):}
\label{Theorem13}
Let $\mathcal{M}$ be an ordinary matroid on the ground set $[n]$. The following are equivalent for a vector $\nu \in \mathbb{R}^n$
\begin{enumerate}
    \item $\nu \in B(\mathcal{M})_{\mathbb{R}}$.
    \item The ordinary matroid $\mathcal{M}_\nu$ does not have loops.
    \item The linear functional $\sum^n_{i=1}x_i\nu_i$ attains its maximum over $P(\mathcal{M})$ in $\partial P(\mathcal{M})\setminus \partial \Delta$. With $\Delta$ being the simplex defined in page $9$.
\end{enumerate}

In the case of a ranked symplectic matroid $\mathcal{S}$, a face of $P(\mathcal{S})$ corresponds to a symplectic matroid but not necessarily to a ranked one. Moreover, it is not true that the symplectic matroid $\mathcal{S}_\omega$ does not have loops for $\omega\in B(\mathcal{S})$. Despite all this, $B(\mathcal{S})$ is still a refinement of the normal fan of $P(\mathcal{S})$.

\subsection{Theorem:}
\label{Theorem14}
Let $\mathcal{S}$ be an ranked symplectic matroid on the ground set $J$ then $B(\mathcal{S})$ is a subfan of the normal fan of $P(\mathcal{S})$.

\subsection*{Proof:}

For every $\omega\in \mathbb{R}^n$, let $\nu\in \mathbb{R}^J$ be the unique non-negative admissible vector with $\nu\in \text{env}^{-1}(\omega)$. By the definition of Bergman fans as order fans of the lattice of flats we have $\omega\in B(\mathcal{S})$ iff $\nu\in  B(\mathcal{M})_\mathbb{R}$.

By Theorem \ref{Theorem1} we have that $\mathcal{M}_\nu$ always contains an admissible basis, and so $\mathcal{M}_\nu$ is an enveloping matroid of $\mathcal{S}_\omega$. Therefore, it is enough to prove that if $\nu_1\in B(\mathcal{M})$ and $\mathcal{S}_{\omega_1}=\mathcal{S}_{\omega_2}$ then $\nu_2\in B(\mathcal{M})$.

Assume $\mathcal{S}_{\omega_1}$ has a loop $a\in J$ as, otherwise, by Theorem \ref{Theorem10} $\mathcal{M}_{\nu_2}$ is also loopless, and $\nu_2\in B(\mathcal{M})$. As $\mathcal{M}_{\nu_1}$ is loopless, it must contain a basis of the form $I\cup \left\{a,a^*\right\}$ and as either $\omega_a=0$ or $\omega_{a^*}=0$ there exist $b\in J$ such that $I\cup \left\{b,a^*\right\},I\cup \left\{b^*,a^*\right\}\in \mathcal{M}_{\nu_1}$. We must then have $I\cup \left\{b,a^*\right\},I\cup \left\{b^*,a^*\right\}\in \mathcal{M}_{\nu_2}$, and thus $\omega_b=\omega_{b^*}=0$ resulting in $I\cup \left\{a,a^*\right\}\in \mathcal{M}_{\nu_2}$. As $\mathcal{M}_{\nu_2}$ is loopless, we get $\nu_2\in B(\mathcal{M})$. \\

The main tool for understanding the algebraic geometric properties of Bergman fans, working over a fixed field
$\mathbb{K}$ is its associated smooth toric variety. For the basics on toric varieties, we refer to \cite{[11]}.

\subsection{Definition:}
\label{Definition12}
To a unimodular fan $\Sigma\subseteq \mathbb{R}^n$ one associates the following smooth toric variety over $\mathbb{K}$:
\[
X(\Sigma)\coloneqq\underset{\sigma\in\Sigma}{\cup}\text{Spec}\mathbb{K}\left[\sigma^\lor \cap \mathbb{Z}\right]
\]
 We observe that the enveloping map defined in section $5$ is induced by the well-defined linear map on the lattices:
 \[
  \text{env}:\mathbb{Z}^J/\left\langle e_{J}\right\rangle \rightarrow \mathbb{Z}^{[n]}
 \]
 \[
  \text{env}(x_1,x_{1^*},...,x_n,x_{n^*})=(x_1-x_{1^*},...,x_n-x_{n^*})
 \]
Moreover we have $\text{env}(\sigma)$ be the cone of $B(\mathcal{S})$ defined by the admissible prefix of its flag of flats for every cone $\sigma$ of $B(\mathcal{M})$, resulting in $\text{env}(B(\mathcal{M}))=B(\mathcal{S})$. Therefore, $\text{env}$ induces a toric morphism between the corresponding toric varieties $\text{env}:X(B(\mathcal{M}))\rightarrow X(B(\mathcal{S}))$.\\

We finish this section with a discussion of the group of Minkowski weights of $B(\mathcal{S})$. Continuing to work with a unimodular fan $\Sigma$ in $N_\mathbb{R}$, let $\Sigma_k$ be the set of k-dimensional cones in $\Sigma$.

\subsection{Definition:}
\label{Definition13}
A function $c:\Sigma_k\rightarrow \mathbb{Z}$ is called a weight of dimension k on $\Sigma$. We say that $c$ is a Minkowski weight of dimension k if for every $k-1$-dimensional cone $\tau$ it satisfies the relation:
\[
\sum_{\tau\subset \sigma}c(\sigma)e_{\sigma/\tau} \text{is contained in the subspace generated by } \tau
\]
With $e_{\sigma/\tau}$ being any representative in $\sigma$ for the generator of the $1$-dimensional lattice $N_\sigma/N_\tau$.

The group of Minkowski weights denoted $MW(\Sigma)$ was first studied in \cite{[12]} with a focus on complete fans. In \cite{[5]} the group of Minkowski weights of dimension k was shown to be isomorphic to $Hom_\mathbb{Z}(A^k(\Sigma),\mathbb{Z})$, the homomorphism group of the Chow cohomology of $\Sigma$. This isomorphism was then used to extend the degree map from complete fans to Bergman fans of ordinary matroids. A key property used was the fact that for a rank $d+1$ ordinary matroid, the group of Minkowski weights of dimension $d$ is isomorphic to $\mathbb{Z}$. We will see that this does not hold for ranked symplectic matroids and explore the structure it does have.

\subsection{Lemma:}
\label{Lemma9}
Let $\mathcal{S}$ be a rank $d+1$ ranked symplectic matroid on the ground set $J$. The constant function $c=1$ is a Minkowski weight of dimension $d$ on $B(\mathcal{S})$.
\subsubsection*{Proof:}
Let $\tau$ be a $d-1$ dimensional cone of $B(\mathcal{S})$. It corresponds to a non-maximal chain of flats $\mathcal{F}=F_1\subsetneq ...\subsetneq F_{d-1}$ of $L(\mathcal{S})\setminus \left\{\emptyset, J\right\}$, missing a flat of rank $r$. Let $\left\{G_1,...G_t\right\}$ be the flats that cover $F_{r-1}$ and, if $r-1\neq d$, be covered by $F_{r}$. Each cone $\sigma$ of dimension $d$ containing $\tau$ contains a unique generating ray not in $\tau$ generated by the lattice point $\sum_{a\in G_i}e^\pm_a$ for some $i\in [t]$. We can therefore choose $e_{\sigma_i/\tau}=\sum_{a\in G_i}e^\pm_a$. If $r=d-1$ then for each $a\in G_i\setminus F_{d-1}$ there exist $j\in [t]$ such that $a^*\in G_j$, thus:
\[
\sum_{\tau\subset \sigma}c(\sigma)e_{\sigma/\tau}=\sum_{i\in [t]}\sum_{a\in G_i}e^\pm_a=t\cdot \sum_{a\in F_{d-1}}e^\pm_a\in \tau
\]
If $r\neq d-1$ then $\left\{G_1\setminus F_{r-1},...G_t\setminus F_{r-1}\right\}$ partition $F_r\setminus F_{r-1}$ and so:
\[
\sum_{\tau\subset \sigma}c(\sigma)e_{\sigma/\tau}=(t-1)\cdot \sum_{a\in F_{r-1}}e^\pm_a+t\cdot \sum_{a\in F_{r}}e^\pm_a\in \tau
\]

\subsection{Corollary:}
In addition to defining a toric variety, the Bergman fan of a ranked symplectic matroid also defines an effective tropical variety with trivial multiplicities.\\

To characterize $MW_d(B(\mathcal{S}))$ for a ranked symplectic matroid $\mathcal{S}$ of rank $d+1$ we investigate two types of pairs of $d$-dimensional cones $\sigma_1,\sigma_2$ of $B(\mathcal{S})$. Let $F_1\subsetneq ...\subsetneq F_{d}$ and $G_1 \subsetneq ...\subsetneq G_d$ be two maximal chains of flats of $L(\mathcal{S})\setminus \left\{\emptyset, J\right\}$ corresponding to $\sigma_1,\sigma_2$. We say that the pair $\sigma_1,\sigma_2$ is of type $1$ if $F_d=G_d$. We say that it is of type $2$ if $F_{d-1}=G_{d-1}$ and $F_d\cap G_d^*\neq \emptyset$. We use the following Lemmas to show that pairs of types $1$ and $2$ determine the structure of $MW_d(B(\mathcal{S}))$.

\subsection{Lemma:}
\label{Lemma10}
Let $\mathcal{S}$ be a rank $d+1$ ranked symplectic matroid on the ground set $J$ and $c$ a Minkowski weight of dimension $d$ on $B(\mathcal{S})$. If $\sigma_1,\sigma_2$ is a pair of type $1$ with $F$ being their shared rank $d$ flat, then $c(\sigma_1)=c(\sigma_2)$.
\subsubsection*{Proof:}
Let $\tau$ be a $d-1$ dimensional cone of $B(\mathcal{S})$ corresponding to a non-maximal chain of flats $\mathcal{F}=F_1\subsetneq ...\subsetneq F_{d-1}=F$ of $L(\mathcal{S})\setminus \left\{\emptyset, J\right\}$ and let $\mathcal{N}$ be the ordinary matroid on the ground set $F$ corresponding to the geometric lattice $[\hat{0},F]\subseteq L(\mathcal{S})$. We observe that the following equality holds:
\[
\sum_{\tau\subset \sigma}c(\sigma)e_{\sigma/\tau}\text{ (mod }\tau)=\sum_{\sigma'}c(\sigma)e_{\sigma'/{\tau'}} \text{ (mod }\tau')
\]
With $\tau'$ being the cone in $B(\mathcal{N})$ corresponding to the non-maximal chain of flats $\mathcal{F}=F_1\subsetneq ...\subsetneq F_{d-2}$ of $L(\mathcal{N})\setminus \left\{\emptyset, F\right\}$ and $\sigma'$ the cone in $B(\mathcal{N})$ corresponding, in the same way, to the generating chain of flats of $\sigma$. As the set of $\sigma'$ is exactly the set of cones containing $\sigma$ we see that the induced function $c':B(\mathcal{N})_{d-1}\rightarrow \mathbb{Z}$ with $c'(\sigma')=c(\sigma)$ is a Minkowski weight of dimension $d-1$ on $B(\mathcal{N})$ and therefore constant, making $c\mid_{\tau\subset \sigma}$ also constant. Using now the fact that $L(\mathcal{S})$ is shellable, we see that $[\hat{0},F]\subseteq L(\mathcal{S})$ is also connected in codimension $1$, implying that $c$ is constant.

\subsection{Lemma:}
\label{Lemma11}
Let $\mathcal{S}$ be a ranked symplectic matroid of rank $d+1$ on the ground set $J$ and $c$ a Minkowski weight of dimension $d$ on $B(\mathcal{S})$. If $\sigma_1,\sigma_2$, corresponding to $F_1\subsetneq ...\subsetneq F_{d}$ and $G_1\subsetneq ...\subsetneq G_d$, is a pair of type $2$ with $F$ being their shared rank $d$ flat then $c(\sigma_1)=c(\sigma_2)$.
\subsubsection*{Proof:}
Using the previous Lemma, we can assume $G_1=F_1,...,G_{d-1}=F_{d-1}$. Let $\tau$ be the cone of dimensions $d$ corresponding to $\mathcal{F}=F_1\subsetneq ...\subsetneq F_{d}=F$ and let $a\in F_d\cap G_d^*$. Observe that in the sum $\sum_{\tau\subset \sigma}c(\sigma)e_{\sigma/\tau}$ the coefficient of $e_{a}$ is exactly $\pm(c(\sigma_1)-c(\sigma_2))$. As $e_{a}$ is orthogonal to the subspace generated by $\tau$ we must have this coefficient be zero.\\

Using what we know so far, we can formulate the following Theorem.

\subsection{Theorem:}
\label{Theorem15}
Let $\mathcal{S}$ be a ranked symplectic matroid of rank $d+1$, $\mathcal{M}$ its enveloping matroid and $\sigma,\sigma'$ be maximal cones in $B(\mathcal{S})$. $c(\sigma)=c(\sigma')$ for all $c\in MW_d(B(\mathcal{S}))$ iff there exists a sequence of $d$-dimensional cones $\sigma=\sigma_1,...,\sigma_k=\sigma'$ such that for each $i\in [k-1]$:
\[
rank_{\mathcal{M}}(F_i\cup F_i^*)\cap (F_{i+1}\cup F_{i+1}^*)=d
\]
with $F_i$ being the maximal flat in the chain of flats corresponding to $\sigma_i$.
\subsubsection*{Proof:}
Let $c\in MW_d(B(\mathcal{S}))$ and let $\sigma,\sigma'$ for which there exists a sequence. It is enough to prove that $c(\sigma)=c(\sigma_2)$, which we do by induction on $1-$ the rank of $F_1\cap F_2$. The basis case $F_1=F_2$ is just a pair of type $1$. Let $a\in F_1\cap F_2^*$ and $a^*\in I$ be an independent set that spans $F_2$, the property above guarantees the existence of such an $a$. Let $\sigma''$ be the maximal cone obtained by replacing $F_2$ with $\overline{I\setminus\left\{a^*\right\}\cup \left\{a\right\}}$ in $\sigma_2$. We find that $\sigma_2,\sigma''$ are a pair of types $2$ and so $c(\sigma_2)=c(\sigma'')$, additionally $c(\sigma)=c(\sigma'')$ by the induction hypotheses. 

On the other hand, if there does not exist such a sequence, then every type $1$ and $2$ pairs $\sigma,\sigma'$ still imply $c(\sigma)=c(\sigma')$, which is enough to make such a $c$ into a $d$-dimensional Minkowski weight.

\subsection{Corollary:}
Let $\mathcal{S}$ be a ranked symplectic matroid of rank $d+1$. If there exists a flat of size $n$ (a transversal in $J$) then $MW_{d}(B(\mathcal{S}))\cong \mathbb{Z}$.\\

In contrast to the case of ordinary matroids, the group of top-dimensional Minkowski weights of a ranked symplectic matroid is not always isomorphic to $\mathbb{Z}$, the most simple example is the free ranked symplectic matroid.

\subsection{Lemma:}
\label{Lemma12}
 $MW_{k-1}(B(U^*_{k,n}))\cong \mathbb{Z}^{n\choose {k-1}}$ 
\subsubsection*{Proof:}
Let $\sigma_1,\sigma_2$ be the maximal cones in $B(U^*_{k,n})$ corresponding to $F_1\subsetneq ...\subsetneq F_{k-1}$ and $G_1\subsetneq ...\subsetneq G_{k-1}$, chains of flats of $L(U^*_{k,n})$. As every base of $\mathcal{M}$, the enveloping matroid, must contain an admissible set of size $k-1$, we have:
\[
rank_{\mathcal{M}}(F_{k-1}\cup F_{k-1}^*)\cap (G_{k-1}\cup G_{k-1}^*)=k \iff G_{k-1}\subsetneq F_{k-1}\cup F_{k-1}^*
\]
Therefore the set of $d$-dimensional Minkowski weights are generated by the set $\left\{A\subseteq [n]\mid  \left|A\right|=k-1\right\}$.

\section{Mason's conjecture for ranked symplectic matroids}

We start by stating Mason's conjecture:

\subsection{Mason's conjecture for ordinary matroids}
For an ordinary matroid on $n$ elements $\mathcal{M}$ with $\mathcal{I}_k$ independent sets of size $k$:
\begin{enumerate}
    \item $\mathcal{I}_k^2\geq \mathcal{I}_{k-1}\cdot \mathcal{I}_{k+1}$
    \item $\mathcal{I}_k^2\geq \left(1+\frac{1}{k} \right)\cdot \mathcal{I}_{k-1}\cdot \mathcal{I}_{k+1}$
    \item $\mathcal{I}_k^2\geq \left(1+\frac{1}{k} \right)\cdot \left(1+\frac{1}{n-k} \right)\cdot \mathcal{I}_{k-1}\cdot \mathcal{I}_{k+1}$
\end{enumerate}

In \cite{[5]} Adiprasito, Huh, and Katz proved (1) log-concavity of the independent set sequence of an ordinary matroid, using techniques from Hodge theory on the Chow ring of a matroid. A proof of (2) was first given in \cite{[14]} by June Huh, Benjamin Schröter, and Botong Wang using the same tools. As we have seen, the Chow ring of a ranked symplectic matroid is much more complicated. An alternative proof, in which (3) was proven, was given independently in \cite{[6]} and \cite{[14]}, both using elementary arithmetic techniques. We will show partial results for ranked symplectic matroids using the two ordinary matroids associated with it. We urge the reader to read both \cite{[6]} and \cite{[14]} for a more complete understanding of the subject and very beautiful theories.\\

Let $\mathcal{S}$ be a ranked symplectic matroid, $\mathcal{M}$ its enveloping matroid, and $\mathcal{N}$ the ordinary matroid of non maximal strongly admissible independent sets of $\mathcal{S}$ defined in Lemma \ref{Lemma4}. For an inadmissible independent set $I$ of $\mathcal{M}$ we have $a\in J$ such that $I'=I\setminus \left\{a\right\}$ is a strongly independent set. We have seen before that for any $b\in I'$ we have that $I'\cup \left\{b^*\right\}$ is independent, and that $(I'\setminus \left\{b\right\})\cup \left\{b^*\right\}\in \mathcal{N}$. Therefore, we can see that under the equivalent relation $I_1\sim I_2$ defined by $I_1\cup I_1^*=I_2\cup I_2^*$ there is a bijection between the independent inadmissible sets of size $k$ of $\mathcal{M}$ and the independent sets of size $k-1$ of $\mathcal{N}$. Moreover, the equivalence classes for independent sets of size $k$ are of uniform size in both matroids. 

\subsection{Lemma:}
\label{Lemma13}
Let $I$ be an inadmissible independent set of $\mathcal{M}$ of size $k$ and $I'$ an independent set of $\mathcal{N}$ of size $k$, then:
\[
\left|[I']_{\sim_\mathcal{N}}\right|=2^{k},\left|[I]_{\sim_\mathcal{M}}\right|=2^{k-1}(k-1)
\]
\subsection*{proof:} 
Since $I'$ is strongly admissible, we have that $(I'\setminus \{ a\} )\cup \{ a^*\}$  is strongly admissible for every $a\in I'$. Hence,
\[
\left|[I']_{\sim_\mathcal{N}}\right|=\left|\{(I'\setminus A)\cup A^*\mid A\subseteq I'\}\right|=2^k
\]
On the other hand, as $I$ is inadmissible, there exists some $a\in J$ such that $I\setminus \{ a\}$  is strongly admissible. Therefore, for any $b\in (I\setminus \{ a\} )$, the set $(I\setminus \{ a\} )\cup \{ b^*\}$ is independent. Moreover, any set in $[I]_{\sim _{\mathcal{M}}}$ must be inadmissible and thus contain a unique pair $a,a^*$. Consequently, we obtain

\[
\left|[I]_{\sim_\mathcal{M}}\right|=\sum_{b\in I\setminus\{a\}}\left|[I\setminus\{a,b\}]_{\sim_\mathcal{N}}\right|=2^{k-1}(k-1)
\]

\subsection{Theorem:}
\label{Theorem16}
Let $\mathcal{S}$ be a ranked symplectic matroid, $\mathcal{M}$ its enveloping matroid and $\mathcal{N}$ the ordinary matroid of its non maximal strongly admissible independent sets. Let $S_k, \mathcal{I}_k, \mathcal{J}_k$ denote the independent sets of size $k$ of $\mathcal{S}$, $\mathcal{M}$ and $\mathcal{N}$, respectively, then:
\[
S_k=\mathcal{I}_k-\frac{k-1}{2}\cdot \mathcal{J}_{k-1}
\]

\subsection*{proof:} 
We already know that $S_k$ is equal to $\mathcal{I}_k$ minus the inadmissible independent sets. By the discussion above, we find that the number of inadmissible independent sets is the same as the number of size $k-1$ independent sets of $\mathcal{N}$ under the equivalent relation. The coefficient $\frac{2}{k-1}$ is a consequence of the size of equivalence classes.\\

We note that from Mason's conjecture for ordinary matroids we have that the sequence of independent sets of a ranked symplectic matroid is a difference of two log concave sequences. we can use this computation to prove log concavity for rank $3$ ranked symplectic matroids:

\subsection{Theorem:}
\label{Theorem17}
The sequence of independent sets for every ranked symplectic matroid of rank $3$ is log-concave. 

\subsection*{proof:} 
The first inequality is obvious and we are left to prove the following:
\[
2n\cdot (\mathcal{I}_3-\mathcal{J}_2)\leq (\mathcal{I}_2-n)^2
\]
We use the following two inequalities: 
\begin{enumerate}
\item 
$2n\cdot \mathcal{I}_3\leq \frac{2}{3}\mathcal{I}_2^2$
\item 
$\mathcal{J}_2\geq n$ 
\end{enumerate}
Where (1) is Mason's conjecture for the enveloping matroid with $k=2$. (2) hold as a rank $3$ matroids has all rank $1$ flats in pairs ($F$ and $F^*$) and has more than one pair. There exist a rank $1$ flat smaller then $\frac{n}{2}$ and any element in it can be extended to more then $n$ strongly independent sets.
Substituting the two inequalities, we obtain:
\[
\frac{1}{3}\mathcal{I}_2^2-2n\mathcal{I}_2+3n^2=(\frac{1}{\sqrt{3}}\mathcal{I}_2-\sqrt{3}n)^2\geq 0    
\]
\[
\Rightarrow (\mathcal{I}_2-n)^2\geq \frac{2}{3}\mathcal{I}_2^2-2n^2 \geq 2n\cdot (\mathcal{I}_3-\mathcal{J}_2)
\]

\address{Einstein Institute of Mathematics, The Hebrew University of Jerusalem, Giv'at Ram, Jerusalem, 91904, Israel}
\email{or.raz1@mail.huji.ac.il}

\newpage
\appendix
\section{Lagrangian orthogonal matroids are ranked}

We show that every Lagrangian orthogonal matroid is a ranked symplectic matroid. We start with some definitions and properties that can be found in \cite{[8]}.

\subsection{Definition:} 
\label{Definition14}
We say an order $\omega$ on $J$ is $D_n$ admissible if and only if it has some admissible set $A$ of $n-1$ elements that are linearly ordered and are the largest elements under $\omega$, $A^*$ are the smallest elements and for $i,j\in A$ such that $i\leq j$ we have $j^*\leq i^*$. The remaining two elements $i,i^*$ are unrelated.
 
An orthogonal matroid on $J$ is a family of equi-numerous admissible subsets of $J$ called bases which contains a unique maximal element with respect to every $D_n$ admissible ordering. An orthogonal matroid of rank $n$ is called a Lagrangian orthogonal matroid.

\subsection{Observation:}
Every orthogonal matroid is a symplectic matroid.

\subsection{Observation:}
All bases of an orthogonal Lagrangian matroid have the same parity in terms of the number of starred elements.

\subsection{Theorem:}
\label{Theorem18}
Every Lagrangian orthogonal matroid is a ranked symplectic matroid. 

\subsection*{Proof:}
Let $\mathcal{B}$ be a Lagrangian orthogonal matroid, by definition it is enough to prove that $\mathcal{M}$, the enveloping ordinary matroid is indeed an ordinary matroid. Let $B_1,B_2\in \mathcal{M}$ be two bases and $a\in B_1$, to prove the basis-exchanged property we consider the following cases:

\begin{enumerate}
\item $B_1,B_2\in \mathcal{B}$, let $b\in B_2\setminus (B_1\setminus \left\{ a\right\})$. We then always have $(B_1\setminus \left\{ a\right\})\cup \left\{ b\right\}\in \mathcal{M}$.
\item $B_1\in \mathcal{B}$ and $B_2\notin \mathcal{B}$, let $b$ be the unique element with $b,b^*\in B_2$ then, as $\left|B_{1}\right|=n$ we must have one of $(B_1\setminus \left\{ a\right\})\cup \left\{ b\right\}$ and $(B_1\setminus \left\{ a\right\})\cup \left\{ b^*\right\}$ be in $\mathcal{M}$.
\item $B_1\notin \mathcal{B}$ and $B_2\in \mathcal{B}$, again let $b$ be the unique element with $b,b^*\in B_1$ and $c\in B_2$ be the unique element with $c,c^*\notin B_1\cup B_1^* $. If $a\in \left\{ b,b^*\right\}$ we have $d\in (B_1\setminus \left\{ a\right\})\setminus B_2^* $, in which case $(B_1\setminus \left\{ a\right\})\cup \left\{ d\right\}\in \mathcal{M}$ or $(B_1\setminus \left\{ a\right\})\cup \left\{ c\right\}=B_2\in \mathcal{M}$. If $a\notin \left\{ b,b^*\right\}$, using the observation that all bases of $\mathcal{B}$ have the same parity, we deduce that we can always add $c$ to $B_1\setminus \left\{ a\right\}$. This is true as we can always extend $(B_1\setminus \left\{ b\right\})$ and $(B_1\setminus \left\{ b^*\right\})$, which are of different parity, to the bases of $\mathcal{B}$.
\item $B_1,B_2\notin \mathcal{B}$, again let $b$ be the unique element with $b,b^*\in B_1$ and $c$ the unique element with $c,c^*\in B_2$. If $a\in \left\{ b,b^*\right\} $ we have one of $(B_1\setminus \left\{ a\right\})\cup \left\{ c\right\}$ and $(B_1\setminus \left\{ a\right\})\cup \left\{ c^*\right\}$ in $\mathcal{M}$. Otherwise, we again have a unique element $d\in B_2$ with $d,d^*\notin B_1\cup B_1^*$ and again we have $(B_1\setminus \left\{ a\right\})\cup \left\{ d\right\}\in \mathcal{M}$ by the same parity argument. 
\end{enumerate}


\begin{thebibliography}{BDP90}


\bibitem {[1]} White, Neil, G-C. Rota, and Neil M. White, eds. Theory of matroids. No. 26. Cambridge University Press, 1986.

\bibitem {[2]} O. Raz, Lattices of flats for symplectic matroids, arXiv:2210.15223 [math.CO]

\bibitem {[3]} Timothy Y. Chow "Symplectic matroids, independent sets, and signed graphs" ,Discrete Mathematics 263 (2003) 35 – 45.

\bibitem {[4]} Andreas Blass, Bruce E. Sagan. "Mobius functions of lattices".  Adv. in Math. 127 (1997), 94-123.

\bibitem {[5]} Karim Adiprasito, June Huh, and Eric Katz. “Hodge theory for combinatorial geometries”. In: Annals of Mathematics 188.2 (2018), pp. 381–452.

\bibitem {[6]} Nima Anari, KuiKui Liu, Shayan Oveis Gharan, Cynthia Vinzant, Log-Concave Polynomials III: Mason’s Ultra-Log-Concavity Conjecture for Independent Sets of Matroids, arXiv:1811.01600.

\bibitem {[7]} Leonid Gurvits. “On multivariate Newton-like inequalities”. In: Advances in Combinatorial Mathematics. Springer, 2009, pp. 61–78.

\bibitem {[8]} Alexandre V. Borovik, Israel M. Gelfand, and Neil White, Coxeter matroid polytopes, Ann. Comb. 1 (1997), no. 2, 123–134, DOI 10.1007/BF02558470. MR1629677 (99d:05016). 

\bibitem {[9]} Feichtner, E.M. and Sturmfels, B. (2005) Matroid Polytopes, Nested Sets and Bergman Fans. Portugaliae Mathematica, 62, 437-468.

\bibitem {[10]} Vladimir Danilov, The geometry of toric varieties, Russian Mathematical Surveys 33 (1978), 97–154.

\bibitem {[11]} William Fulton, Introduction to Toric Varieties, Annals of Mathematics Studies, 131, Princeton University Press,
Princeton, NJ, 1993.

\bibitem {[12]} William Fulton and Bernd Sturmfels, Intersection theory on toric varieties, Topology 36 (1997), no. 2, 335–353

\bibitem {[13]} Federico Ardila and Caroline Klivans, The Bergman complex of a matroid and phylogenetic trees, Journal of Combinatorial Theory Series B 96 (2006), no. 1, 38–49.

\bibitem {[14]} June Huh, Benjamin Schröter, Botong Wang, Correlation bounds for fields and matroids. J. Eur. Math. Soc. 24 (2022), no. 4, pp. 1335–1351.

\bibitem {[15]} Brändén Petter and June Huh. 2020. Lorentzian polynomials. Annals of Mathematics 192(3):821–891.

\bibitem {[16]} A. Bjorner, A. M. Garsia and R.P. Stanley, An introduction to Cohen-Macaulay partially ordered sets, in " Ordered Sets" (ed. I. Rival), Reidel, Dordrecht, 1982, pp. 583-615.

\bibitem {[17]} Chris Godsil, An Introduction to the Moebius Function, arXiv:1803.06664 [math.CO]

\end{thebibliography}
\end{document}